\documentclass{amsart}
\usepackage{amsmath, amsthm, amssymb}

\newtheorem{theorem}{Theorem}[section]
\newtheorem{lemma}[theorem]{Lemma}
\newtheorem{proposition}[theorem]{Proposition}
\newtheorem{corollary}[theorem]{Corollary}

\newtheorem{remark}[theorem]{Remark}

\newcommand{\bfind}[1]{\index{#1}{\bf #1}}
\newcommand{\n}{\par\noindent}

\newcommand{\sn}{\par\smallskip\noindent}
\newcommand{\mn}{\par\medskip\noindent}

\newcommand{\pars}{\par\smallskip}
\newcommand{\parm}{\par\medskip}

\newcommand{\cal}{\mathcal}
\newcommand{\cG}{{\mathcal G}}

\newcommand{\chara}{\mbox{\rm char}\,}
\newcommand{\trdeg}{\mbox{\rm trdeg}\,}

\newcommand{\subsetuneq}{\mathrel{\raisebox{.8ex}{\footnotesize%
$\displaystyle\mathop{\subset}_{\not=}$}}}

\newcommand{\pH}{\mathop{\raisebox{-.2ex}{\rule[-1pt]{0pt}{1pt}%
\mbox{\large\bf H}}}}

\newcommand{\supp}{\mbox{\rm supp}\,}
\newcommand{\N}{\mathbb N}
\newcommand{\Q}{\mathbb Q}
\newcommand{\R}{\mathbb R}
\newcommand{\Z}{\mathbb Z}

\newcommand{\F}{\mathbb F}

\begin{document}

\title[Methods for the classification of cuts]{Selected methods for the classification of cuts, and their applications}
\author{F.-V.~Kuhlmann}
\address{Institute of Mathematics, University of Silesia at Katowice, Bankowa 14,
40-007 Katowice, Poland}
\email{fvk@math.us.edu.pl}

\thanks{The author would like to thank Katarzyna Kuhlmann and Marcus Tressl for useful suggestions.}


\date{28.\ 2.\ 2017}

\begin{abstract}
We consider four approaches to the analysis of cuts in ordered abelian groups and ordered fields, their
interconnection, and various applications. The notions we discuss are: ball cuts, invariance group, invariance
valuation ring, and cut cofinality.
\end{abstract}

\maketitle

%
%
\section{Introduction}
In these notes we deal with (Dedekind) cuts in ordered abelian groups and in ordered fields. (For the definition
of the notion of a cut and other notions used in this Introduction, see Section~\ref{sectprel}.)
We introduce the reader to four approaches to their classification, the links between them, and several
applications. The reader should observe that a cut in an ordered field is at the same time a cut in its additive
group. Hence even in the case that one is predominantly interested in cuts in ordered fields, up to a certain point
their study can be fruitfully carried out in the setting of ordered abelian groups and does not need to make use of
the field (or ring) multiplication. At the same time the reader should keep in mind that ordered abelian groups
appear in field theory also as the value groups of valuations. In this case, cuts in the value group can for
instance be
generated by pseudo Cauchy sequences in the valued field. If the field is ordered and its valuation is the natural
valuation induced by the ordering (see below), then it is essential to study the connection between cuts in the
field and induced cuts in the value group.

\pars
The first approach to the classification of cuts is to ask whether a cut in an ordered abelian group is the upper
or lower edge of a convex subgroup, or of a coset thereof. This has been used and studied more or less
explicitly by many authors, and various names have been given to such cuts. We call them \bfind{ball cuts}.
They appear implicitly or explicitly, sometimes with surprisingly different
definitions, in \cite{G1,G2,G3,GP,Pe} for the study of cuts in ordered fields, in \cite{[W]} for the study of
cuts in ordered abelian groups, and in several other papers cited in the references.
Ball cuts will be introduced and discussed in Section~\ref{sectballcuts}.

\pars
Spaces of $\R$-places (i.e., places with residue fields embeddable in $\R$) of ordered fields are not well
understood. It is a longstanding open problem which topological spaces appear as spaces of $\R$-places. Recently,
ball cuts have been used to study these spaces (cf.\ \cite{KK1,[KMO],Kasia,KKo}). We describe some results in
Section~\ref{sectRpla}.

\pars
Two well known deep open problems in positive characteristic are:
\sn
1) resolution of singularities in arbitrary dimension,
\sn
2) decidability of the field $\F_p((t))$ of Laurent series over a finite field.
\sn
Both problems are connected with the structure theory of valued function fields of positive characteristic $p$. The
main obstruction here is the phenomenon of the \bfind{defect}. Via ramification theory, the study of it can be
reduced to the study of purely inseparable extensions and of Galois extensions of degree $p$. Ball cuts are
essential for a classification of Galois extensions with nontrivial defect which is introduced in \cite{[Ku4]}.
It will be discussed in Section~\ref{sectASe}.

\parm
How ``broad'' is a given cut? One way to answer this question is to associate to the cut the maximal set of
elements that can be added to the cut sets without changing the cut. This set turns out to be a convex subgroup
of the ordered group; we call it the \bfind{(additive) invariance group}. This notion was introduced by the author
in his thesis (\cite{[Ku1]}) in order to handle valued fields with non-archimedean ordered value groups in
connection with the model theory of valued fields. Invariance groups were also introduced by M.~Tressl in his
thesis (\cite{[T]}), this time for the study of the model theory of ordered fields (a quick example for Tressl's use
of invariance groups is given in Section~\ref{sectTres}). Later, A.~Fornasiero and
M.~Mamino (\cite{F}) used them in a detailed investigation of cuts of ordered abelian groups, which they then
applied to study so-called double ordered monoids. Moreover, they have been implicitly used by several other
authors (e.g.\ by R.~Rolland in \cite{[Ro]}), or even explicitly defined under different names (e.g.\ by
F.~Wehrung in \cite{[W]}, and by D.~Kijima and M.~Nishi in \cite{[KN]}). A main link
to ball cuts is the fact that invariance groups can help to identify them (see Theorem~\ref{ball-invg} below).
We will discuss invariance groups in Section~\ref{sectinvgp}.

\pars
In mathematics, objects that are maximal with respect to a certain property are often of particular interest. In
valuation theory, this is so for \bfind{maximal valued fields}, which have the property that every proper extension
will necessarily enlarge value group or residue field. Ordered abelian groups and fields carry \bfind{natural
valuations} which are canonically derived from their ordering. In Section~\ref{sectmaxf} we describe a
characterization of certain ordered fields maximal with respect to their natural valuation, given by
Kijima and Nishi; it makes essential use
of invariance groups. Further, we discuss a generalization of this result, due to H.-J.~H\"uper, to the case of
valuations whose valuation rings are convex under the given ordering.

In Section~\ref{sectmaxg} we discuss an analogue for the case of ordered abelian groups. We present the Cohen-Goffman Theorem and a related result by P.~Ehrlich, both of which (implicitly) use invariance groups.

\parm
The author's attention was drawn to the importance of invariance groups in the study of ordered fields by a
question of J.~Madden. During the Special Semester in Real Algebraic Geometry and Ordered
Structures, Baton Rouge 1996, Madden showed him the definition of what we now
call the \bfind{invariance valuation ring} and asked for the meaning of it.
The answer to his question was first given in the manuscript \cite{[Ku3]}. Again,
the invariance valuation ring was independently introduced and applied by M.~Tressl (see \cite{[T],[T1],[T2]}).
The construction of invariance valuation rings appears already in Rolland's paper \cite{[Ro]},
but not in full generality.

\pars
Looking at a cut in an ordered field $K$, one may ask whether it
originates in some way from a cut in the residue field of $K$ with
respect to some real place. That is, one would like to know  whether the
cut can be translated into some ``normal position'' such that for some
convex valuation ring ${\cal O}$ of $K$ with maximal ideal
${\cal M}$, it induces a (Dedekind) cut in the ordered residue field ${\cal O}/{\cal M}$
via the residue map. If so, one would like to determine how this translation can be done.
The invariance valuation ring is a key tool to answer these questions (see
Section~\ref{sectproj}).

Further, in the paper \cite{JSW} F.~Jahnke, P.~Simon and E.~Walsberg introduce certain invariance valuation rings
to exhibit definable valuations in ordered fields which are not dense in their real closure. The details will be
discussed in Section~\ref{sectdefval}.

\parm
A remark by M.~Marshall made it clear to the author of these notes that some of his earlier results were
actually a special case of a more general setting which we will now sketch. Every ordered field $K$ has a
\bfind{natural valuation} $v$, whose residue field
is an archimedean ordered field; its valuation ring ${\cal O}_v$ is the smallest of all convex valuation rings of
$K$. Then every convex subgroup of the ordered additive group of $K$ is an ${\cal O}_v$-module.

In the present paper, we use the additive Krull notation for valuations, i.e., the ultrametric triangle law reads
$v(a+b)\geq \min\{v(a),v(b)\}$ and the value group is an additively written ordered abelian group whose
nonnegative elements are precisely the values of the nonzero elements of ${\cal O}_v\,$. In this notation,
the map
\[
M\>\mapsto\>(vK\setminus vM,vM)\;,
\]
where $vK$ is the value group of $(K,v)$ and $vM:=\{va\mid 0\ne a\in
M\}$, is a bijection between the convex subgroups $M$ of $K$ and
the cuts in the value group $vK$. This holds more generally for any
(Krull) valuation $v$ of an arbitrary field $K$ and the set of all
${\cal O}_v$-modules $M\subseteq K$. Information about $M$ can be read off
from the invariance group of the cut $(vK\setminus vM,vM)$. One can also
define the invariance valuation ring of an ${\cal O}_v$-module. The
invariance valuation ring of a cut can then be understood as the invariance valuation ring
of the invariance group of the cut.

\pars
Tressl introduced the author to the definition and main properties of the
\bfind{multiplicative invariance group} of a cut in an ordered field, that is, the invariance group of the cut
taken in the multiplicative group of the field. For its properties, see \cite{[Knew]}, where
a detailed study of ball cuts, invariance groups and invariance valuation rings is presented. Detailed studies of
cuts using these concepts appear also in Tressl's papers \cite{[T],[T1],[T2]} and T.~G\"uldenberg's thesis \cite{Gu}.

\parm
After ball cuts, invariance group and invariance valuation ring, the fourth approach to the study of Dedekind cuts
is to consider the pair of cardinal numbers
$(\kappa,\lambda)$ where $\kappa$ is the cofinality of the lower cut set and $\lambda$ is the coinitiality of
the upper cut set. Recall that the coinitiality of a linearly ordered set is
the cofinality of this set under the reversed ordering. Recall further
that cofinalities and coinitialities of ordered sets are regular cardinals.
We call $(\kappa,\lambda)$ the \bfind{cofinality} of the cut; also the name \bfind{character} has been used in
the literature.

In his groundbreaking and comprehensive work, Hausdorff constructs for any given collection of cofinalities
$(\kappa,\lambda)$ which satisfy some necessary conditions, a totally ordered set where this collection is exactly
the set of cofinalities of the cuts appearing in this ordering. One aim of the already cited paper \cite{[Ro]} of
Rolland is to construct ordered fields which realize a prescribed set of cut cofinalities.

\pars
A much studied property of ordered abelian groups or fields is that of being an $\eta_\alpha$-set, which is
equivalent to the absence of cuts of cofinality $(\kappa,\lambda)$ with both$\kappa$ and $\lambda$ smaller than
$\aleph_\alpha\,$. We discuss a characterization of such ordered abelian groups and fields, due to N.~Alling, in
Section~\ref{sectea}.

In Section~\ref{sectgp}, we present some work of N.~Yu.~Galanova and G.~G.~Pestov which involves cut cofinalities
and ball cuts.

\pars
More recently, a new aspect of cut cofinalities has been discovered. Transferring the concept of spherical
completeness from ultrametric spaces to other spaces equipped with distances or topologies, the authors of
\cite{KK2} asked the question whether there are ordered fields, apart from the reals themselves, in which every
chain of closed bounded intervals has a nonempty intersection. This happens exactly when all appearing cut
cofinalities $(\kappa,\lambda)$ satisfy $\kappa\ne\lambda$. The positive answer to the question was first given by
S.~Shelah in the paper \cite{[S]}. In joint work with Shelah and K.~Kuhlmann, the author of these notes gave an
alternative construction and a complete characterization of such fields in \cite{KKS}.
We will discuss some details in Section~\ref{sectsc}.

\parm
These notes are not intended to be a comprehensive survey on the general theory of cuts. However, the author hopes
that they will initiate discussion and feedback so that more comprehensive information can be gathered and later be
put together in a monograph on cuts.

\mn
%
%
\section{Notation and preliminaries}                   \label{sectprel}
For general background from valuation theory, we recommend \cite{[EP]}. For background on ordered fields, see
\cite{L,Pr}.

\sn
\subsection{Cuts}
Take any ordered set $(S,<)\,$ ({\it by ``ordered'', we will always
mean ``totally ordered''}). If $S_1, S_2$ are nonempty subsets of $S$ and
$a\in S$, we will write $a<S_2$ if $a<b$ for all $b\in S_2$, and
we will write $S_1<S_2$ if $a<S_2$ for all $a\in S_1$.

A subset $S'$ of $S$ is called \bfind{convex} {\bf in} $(S,<)$ if for
every two elements $a,b\in S'$ and every $c\in S$ such that $a\leq c\leq
b$, it follows that $c\in S'$.
A subset $S_1$ of $S$ is an \bfind{initial segment} {\bf of} $S$ if
for every $a\in S_1$ and every $c\in S$ with $c \leq a$, it follows
that $c\in S_1\,$. Symmetrically, $S_2$ is a \bfind{final segment}
{\bf of} $S$ if for every $a \in S_2$ and every $c\in S$ with $c\geq
a$, it follows that $c\in S_2\,$. Note that $S_1$ is an initial
segment of $S$ if and only if $S_1$ is convex and $S_1< S
\setminus S_1\,$. Note also that $\emptyset<S$ and $S<\emptyset$ by
definition; so $\emptyset$ is an initial segment as well as a final
segment of $S$.

If $S_1\subseteq S$ and $S_2\subseteq S$ are such that $S_1< S_2$ and $S=S_1\cup S_2$, then we will call
$(S_1,S_2)$ a \bfind{cut} in $S$. Then $S_1$ is an initial segment of $S$, $S_2$ is a final segment of $S$,
and the intersection of $S_1$ and $S_2$ is empty. We write $\Lambda^L=S_1\,$,
$\Lambda^R=S_2\,$, and
\[
\Lambda\>=\>(\Lambda^L,\Lambda^R)\;.
\]
A cut $(\Lambda^L,\Lambda^R)$ with $\Lambda^L\ne\emptyset$ and
$\Lambda^R\ne\emptyset$ is called a \bfind{Dedekind cut}. If $\Lambda$ is a cut in $S$,
$(T,<)$ is an extension of $(S,<)$ and $a \in T$ is such that
$\Lambda^L\leq a\leq\Lambda^R$, then we will say that \bfind{$a$ realizes $\Lambda$} (in $T$).

For any subset $M\subseteq S$, we let $M^+$ denote the cut
\[
M^+\>=\>(\{s\in S\mid \exists m\in M: s\leq m\}\,,\,\{s\in S\mid s>M\})\;.
\]
That is, if $M^+=(\Lambda^L,\Lambda^R)$ then $\Lambda^L$ is the least
initial segment of $S$ which contains $M$, and $\Lambda^R$ is the
largest final segment having empty intersection with $M$. If $M=
\emptyset$ then $\Lambda^L=\emptyset$ and $\Lambda^R=M$, and if $M=S$,
then $\Lambda^L=M$ and $\Lambda^R=\emptyset$. Symmetrically, we set
\[
M^-\>=\>(\{s\in S\mid s<M\}\,,\,
\{s\in S\mid \exists m\in M: s\geq m\})\;.
\]
That is, if $M^-=(\Lambda^L,\Lambda^R)$ then $\Lambda^L$ is the largest
initial segment having empty intersection with $M$, and $\Lambda^R$ is the least
final segment of $S$ which contains $M$. If $M=\emptyset$ then
$\Lambda^L=M$ and $\Lambda^R=\emptyset$, and if $M=S$, then
$\Lambda^L=\emptyset$ and $\Lambda^R=M$.

If $M=\{a\}$, we will write $a^+$ instead of $\{a\}^+$ and $a^-$ instead of $\{a\}^-$. These two cuts are called
\bfind{principal}. Hence if $M$ has a largest element $a$, then $M^+=a^+$ is principal, and if $M$ has a smallest
element $a$, then $M^-=a^-$ is principal. The cut $(\Lambda^L,\Lambda^R)$ is principal if and only if
$\Lambda^L$ has a largest element or $\Lambda^R$ has a smallest element. In the literature, a principal cut is
also called \bfind{realized} or \bfind{filled}, a non-principal cut is called a \bfind{gap}, and a cut for which
$\Lambda^L$ has a largest element {\it and} $\Lambda^R$ has a smallest element is called a \bfind{jump}. In
\cite{Eh}, Ehrlich calls a cut $(\Lambda^L,\Lambda^R)$ \bfind{continuous} if $\Lambda^L$ is principal but not a jump.

\sn
\subsection{Valuation theory}   \label{sectvalth}
Take an ordered abelian group $G$. Two elements $a,b$ are \bfind{archimedean equivalent} if there is some
$n\in\N$ such that $n|a|\geq |b|$ and $n|b|\geq |a|$. The equivalence
class of $a$ is called \bfind{archimedean class of} $a$ and is denoted
by $[a]$. The set $\{[a]\mid 0\ne a\in G\}$ is totally ordered by
setting $[a]< [b]$ if and only if $|a|>n|b|$ for all $n\in\N$. Then the class of $0$ is the largest element in the
set, and it only contains
the element $0$. The map $v:\>a\mapsto va:=[a]$ is the \bfind{natural valuation} of $G$. It satisfies
the triangle inequality $v(a+b)\geq \min\{va,vb\}$ and $v(-a)=va$. We call $\{va\mid 0\ne g\in G\}$ the
\bfind{value set of $G$ (under $v$)}.

If $G$ is the additive group of an ordered
field $K$, then by setting $[a]+ [b]:= [ab]$ we obtain an addition on the set of archimedean classes that is
compatible with the ordering, and the natural valuation becomes a field (Krull) valuation.

\parm
Take any extension $(L|K,v)$ of valued fields, that is, an extension
$L|K$ of fields and a valuation $v$ on $L$. By $vL$ and $vK$ we denote
the value groups of $v$ on $L$ and on $K$, and by $Lv$ and $Kv$ the
residue fields of $v$ on $L$ and on $K$, respectively. Similarly, $vz$
and $zv$ denote the value and the residue of an element $z$ under $v$.

\pars
A valued field $(K,v)$ is called \bfind{henselian} if the extension of $v$ to every algebraic extension field $L$
of $K$ is unique, or equivalently, $(K,v)$ satisfies Hensel's Lemma. A \bfind{henselization} of $(K,v)$ is an
algebraic extension of $(K,v)$ which is henselian and can be embedded over $K$ in every other henselian extension
field of $(K,v)$. Henselizations exist and are unique up to valuation preserving isomorphism over $K$. Therefore,
we will speak of {\it the} henselization of $(K,v)$ and denote it by $K^h$.

\pars
Assume that $L|K$ is finite and the extension of $v$ from $K$ to $L$ is
unique. Then the Lemma of Ostrowski says that
\begin{equation}                            \label{LoO}
[L:K]\;=\; p^\nu \cdot (vL:vK)\cdot [Lv:Kv] \;\;\;\mbox{ with }\nu\geq 0
\end{equation}
where $p$ is the \textbf{characteristic exponent} of $Kv$, that is,
$p=\chara Kv$ if this is positive, and $p=1$ otherwise.
The factor ${\rm d}=p^\nu$ is called the \bfind{defect} of the
extension $(L|K,v)$. If d$\>=1$, then we call $(L|K,v)$ a
\bfind{defectless extension}; otherwise, we call it a \textbf{defect
extension}. Note that $(L|K,v)$ is always defectless if $\chara Kv=0$.

We call a henselian field $(K,v)$ a \bfind{defectless field} if every finite extension of $(K,v)$ is defectless.
An arbitrary field is called a defectless field if its henselization is defectless.

\pars
The extension $(L|K,v)$ is \bfind{immediate} if for each $z\in L\setminus K$ there is $c\in K$ such that
$v(z-c)>vz$; this holds if and only if the canonical embeddings of $vK$ in $vL$ and of $Kv$ in $Lv$ are onto.

\pars
For $z\in L$, we define
\[
v(z-K)\>:=\> \{v(z-c)\mid c\in K\}\>\subseteq\> vL\cup\{\infty\}\;.
\]
If $(L|K,v)$ is immediate, then $v(z-K)$ is a subset of $vK$ without a maximal element, and even more, it is an
initial segment.

Immediate extensions of valued abelian groups can be defined as in the case of valued fields. Valued abelian
groups and valued fields are called \bfind{maximal} if they do not admit proper immediate extensions.

\sn
\subsection{Pseudo Cauchy sequences}             \label{sectpCs}
A \bfind{pseudo Cauchy sequence} in a valued abelian group or field is a sequence $(a_\nu)_{\nu<\lambda}$ of
elements, indexed by a limit ordinal $\lambda$ (which is called the \bfind{length} of the sequence), such that
for all $\rho<\sigma<\tau<\lambda$,
\[
v(a_\sigma-a_\rho)\><\>v(a_\tau-a_\sigma)\>.
\]
In this case, $v(a_\sigma-a_\rho)=v(a_{\rho+1}-a_\rho)\,$. If $(a_\nu)_{\nu<\lambda}$ is a pseudo Cauchy sequence,
then the sequence of values $(v(a_{\nu+1}-a_\nu))_{\nu<\lambda}$ is strictly increasing. The set $\{b\mid
\forall\nu<\lambda:\> v(b)>v(a_{\nu+1}-a_\nu)\}$ is called the \bfind{breadth} of the sequence
$(a_\nu)_{\nu<\lambda}$.
An element $a$ (in some valued extension group or field) is a \bfind{limit} of the sequence if $v(a-a_\nu)=
v(a_{\nu+1}-a_\nu)$ for all $\nu<\lambda$. If $a$ is a limit of the sequence, than also $a'$ is a limit if and
only if $a-a'$ is an element of the breadth.

A valued abelian group or field is called \bfind{spherically complete} if it admits a limit for every pseudo Cauchy
sequence. If the valued abelian group $G'$ is an immediate extension of the valued group $G$, then every element
$a\in G' \setminus G$ is the limit of a pseudo Cauchy sequence in $G$ that does not have a limit in $G$. Hence,
every spherically complete valued abelian group or field is maximal. The converse is also true; in the case of
valued fields this is shown by I.~Kaplansky in \cite{[Ka]}, where the theory of pseudo Cauchy sequences (which he
calls ``pseudo-convergent sequences'') is nicely laid out.

If $\alpha$ is an ordinal, then $G$ is called \bfind{$\alpha$-maximal} if every pseudo Cauchy sequence in $G$ of
length less than $\aleph_\alpha$ has a limit in $G$.

%
%
\subsection{Hahn products and power series fields}                  \label{secthahn}
Given a linearly ordered index set $I$ and for every $\gamma\in I$ an
arbitrary abelian group $C_\gamma\,$, we define a group called the \bfind{Hahn product} (also called
\bfind{Hahn group}), denoted by $\pH_{\gamma\in I} C_\gamma\,$.
Consider the product $\prod_{\gamma\in I} C_\gamma$ and an element $c=
(c_\gamma)_{\gamma \in I}$ of this group. Then the \bfind{support} of
$c$ is the set $\supp c:=\{\gamma\in I\mid c_\gamma \not=0\}$. As a
set, the Hahn product is the subset of $\prod_{\gamma\in I} C_\gamma$
containing all elements whose support is a wellordered subset of $I$,
that is, every nonempty subset of the support has a minimal element. The Hahn product is
a subgroup of the product group. Indeed, the support of the (componentwise) sum of two
elements is contained in the union of their supports, and the union of
two wellordered sets is again wellordered.

The support of every nonzero element $c$ in the Hahn product
has a minimal element $\gamma_0\,$. This enables us to define a group
valuation by setting $vc=\gamma_0$ and $v0=\infty$; this is called the \bfind{canonical valuation}
of the Hahn product $\pH_{\gamma\in I} C_\gamma\,$.

\pars
If the $C_\gamma$ are (not necessarily archimedean) ordered
abelian groups, we obtain the \bfind{ordered Hahn product}, also called
\bfind{lexicographic product}, where the ordering is defined as
follows. Given a nonzero element $c= (c_\gamma)_{\gamma \in I}\,$, let
$\gamma_0$ be the minimal element of its support. Then we take $c>0$ if
and only if $c_{\gamma_0}>0$. If all $C_\gamma$ are archimedean
ordered, then the canonical valuation of the Hahn product coincides with the
natural valuation of the ordered Hahn product. The \bfind{Hahn Embedding Theorem} states that
every ordered abelian group $G$ can be embedded in the Hahn product with its set of
archimedean classes as index sets and all $C_\gamma$ equal to the ordered group of real numbers.

Take any ordered abelian group $G$. If $H\subsetuneq H'$ are convex subgroups of $G$ such that the ordering induced
on $H'/H$ is archimedean (and hence $H'/H$ can be seen as an ordered subgroup of the reals), then $H'/H$ is called
an \bfind{archimedean component of $G$}. If $G=\pH_{\gamma\in I} C_\gamma$ and all $C_\gamma$ are archimedean
ordered, then the $C_\gamma$ are precisely the archimedean components of $G$.

\pars
Take a field $k$ and an ordered abelian group $G$. Then $k((G)):=\pH_{\gamma\in G} k$ is a valued abelian group.
Since all supports are wellordered, a multiplication can be defined as follows: $(c_g)_{g\in G} \cdot
(c'_g)_{g\in G} = (\sum_{h+h'=g} c_h\cdot c'_{h'})_{g\in G}\,$. Then $k((G))$ becomes a valued field, called a
\bfind{power series field}. The canonical valuation of the underlying Hahn product makes it a valued
field with value group $G$ and residue field $k$.

Under their canonical valuation, all Hahn products and all power series fields are spherically complete
and hence maximal. All maximal fields with residue fields of
characteristic $0$ are power series fields, but for positive residue characteristic this is not true.

\mn
%
%
\section{Ball cuts}                   \label{sectballcuts}
We say that a cut $\Lambda=(\Lambda^L,\Lambda^R)$ in an ordered abelian group is a \bfind{group$\,^+$-cut} if it is
induced by the upper edge of a convex subgroup $H$ of $G$, i.e., if $\Lambda=H^+$. We will say that $\Lambda$ is a
\bfind{group$\,^-$-cut} if it is induced by the lower edge of a convex
subgroup $H$ of $G$, i.e., if $\Lambda=H^-$. In both cases, we will call $\Lambda$ a \bfind{group-cut}.
Note that $0^+$ and $0^-$ are the only principal group-cuts.
We call $\Lambda$ a \bfind{ball$\,^+$-cut} (or a \bfind{ball$\,^-$-cut}) if it is induced by the upper edge (or
lower edge, respectively) of some coset of a convex subgroup $H$ of $G$, i.e., if it is of the form $(g+H)^+$ (or
$(g+H)^-$, respectively) for some $g\in G$. Ball$\,^+$-cuts and ball$\,^-$-cuts are called \bfind{ball-cuts}, and
cosets $H+g$ of convex subgroups $H$ are also called \bfind{balls}. Note that all group-cuts are ball-cuts.

\pars
Ball cuts are called \bfind{asymmetric cuts} in \cite{G1,G2,G3,GP,Pe}. This name is unfortunate; it may have been
chosen by the authors after they observed that there are no cuts in ordered fields that are at the same time a
ball$^-$-cut and a ball$^+$-cut. But the situation is different in ordered abelian groups, as the following example
shows. Consider the lexicographic ordering on $\Z\times\Z$. Then $\{(0,m)\mid m\in\Z\}^+=\{(1,m)\mid m\in\Z\}^-$.

\pars
In Tressl's paper \cite{[T2]}, the ball$^+$-cuts are the cuts with \bfind{signature} $1$, and the ball$^-$-cuts
are the cuts with signature $-1$. All non-ball cuts have signature $0$. G\"uldenberg also uses signatures in his
thesis \cite{Gu}, but defines them in a slightly different way.

\sn
\subsection{Monoids of cuts}                  
On the set of cuts in an ordered abelian group, an addition can be defined in various ways. The two immediately
obvious ways to define $\Lambda_1+\Lambda_2$ are the following:
\sn
1) \ set $\Lambda_1+\Lambda_2:=(\Lambda_1^L+\Lambda_2^L)^+=\{\alpha+\beta\mid \alpha\in \Lambda_1^L, \beta\in \Lambda_2^L\}^+$,
\n
2) \ set $\Lambda_1+\Lambda_2:=(\Lambda_1^R+\Lambda_2^R)^-$.
\sn
The two additions are usually not the same, but their properties are very similar. The following fact is easy to
prove:
\sn
{\it The idempotent elements in these monoids are precisely the group cuts.}

\pars
Monoids of cuts are studied in \cite{[W]}, \cite{Gu}, \cite{F} and \cite{FKK}. In the latter paper, the results are
used for the intrinsic construction (without the use of embeddings in power series fields) of towers of complements
to all (possibly fractional) ideals of the valuation ring in henselian valued fields of residue characteristic 0,
and in Kaplansky fields (i.e., valued fields satisfying ``Hypothesis A'' in \cite{[Ka]}) which do not admit proper
immediate algebraic extensions. They are also used by N.~Alling in \cite{A} for the characterization of
$\eta_\alpha$ ordered abelian groups and fields (see Section~\ref{sectea} below). Alling gives credit to
A.~H.~Clifford (\cite{C}) for introducing the monoid structure (but it had very probably already been observed
before, when Dedekind completions of ordered abelian groups were considered).

\sn
\subsection{Applications to spaces of $\R$-places}                  \label{sectRpla}
For any formally real (i.e., orderable) field $K$, the question arises which orderings induce the same natural
valuations. The places
associated with natural valuations are called $\R$-places as their residue fields are archimedean ordered and can
thus be embedded in $\R$. We will therefore always assume that the residue field of an $\R$-place is a
subfield of $\R$. The above question was answered in \cite{[KMO]} for an interesting special case.

Take a real closed field $R$. There is a one-to-one correspondence between orderings $P$ of $R(X)$ and
cuts of $R$ (see \cite{g}). The cut $\Lambda_P=(\Lambda_P^L,\Lambda_P^R)$ corresponding
to $P$ is given by $\Lambda_P^L =\{a \in R\mid a <_P X\}$ and $\Lambda_P^R =\{b \in R
\mid b >_P X\}$. Conversely, if $\Lambda$ is a cut in $R$, then the set
\[
P = \{ f \in R(X)\mid  \exists a\in \Lambda^L \;\exists b \in \Lambda^R\; \forall c \in
(a,b) :\; f(c) \in \dot R^2\}
\]
is an ordering of $R(X)$ with $\Lambda_P = \Lambda$ (here, $\dot R=R\setminus\{0\}$).

\pars
In \cite{[KMO]} the following result is proved:

\begin{theorem}                  \label{2to1}
Two distinct orderings of $R(X)$ induce the same $\R$-place if and only if they correspond to the upper and lower
edges of the same ball, that is, there is a convex subgroup $H$ of the additive group of $R$ and $c\in R$ such
that one of the places corresponds to $(H+c)^-$ and the other to $(H+c)^+$.
\end{theorem}

This means that the space of $\R$-places of $R(X)$ is obtained from the line $R$ by identifying the upper and lower
edges of balls. If this is done for $R=\R$ then we obtain the circle (up to homeomorphism). But if $R$ is a
non-archimedean ordered field, then the structure is much more complex; it may be thought of as an infinite pearl
necklace in which every pearl contains a pearl necklace that is similar to the whole necklace. The rich
self-similarities of this space have been exhibited in \cite{Kasia} by observing that the transformations
$a\mapsto a+c$, $a\mapsto ca$ and $a\mapsto a^{-1}$ all transform balls into balls.

\pars
The following result is also proved in \cite{[KMO]}:

\begin{theorem}
Take an ordering on $R(X)$ which extends the ordering of $R$, and take $v$ to be the natural valuation on $R(X)$
w.r.t.\ this ordering. Then $X$ induces in $R$ a cut of the form $(c+H)^-$ or $(c+H)^+$ (as in Theorem~\ref{2to1})
if and only if $vR\subsetuneq vR(X)$, and if the former is the case, then $v(X-c)$ is rationally independent over
$vR$.
\end{theorem}

From this theorem we conclude that a cut of $R$ is a ball cut if and only if the natural value group $vR(X)$ of
the corresponding ordering on $R(X)$ satisfies $[vR(X):2vR(X)] = 2$.

\pars
In the paper \cite{KKo} P.~Koprowski and K.~Kuhlmann consider the more general case of an algebraic function
field $F$ of transcendence degree 1 over a
real closed field $R$. Choose any smooth projective model of $F$, i.e., a smooth, projective
algebraic curve over $R$ with function field $F$. In \cite{[Kn]} M.~Knebusch shows that the curve consists of
finitely many semialgebraic connected components, each of which can be endowed with a cyclic order. In \cite{KKo}
this is used to define cuts in these components; the collection of all of them is taken to be the set of cuts on the
curve. The following result is proved:

\begin{theorem}                            \label{homeo}
The space of all cuts on the curve (endowed with the order topology) is homeomorphic to the space of all orderings
on $F$ (endowed with the Harrison topology).
\end{theorem}

Take any ordering of $F$ and let $v$ denote the natural valuation of $F$ w.r.t.\ this ordering. Note that the value
group $vR$ of $v$ on $R$ is divisible since $R$ is real closed. Therefore, as $\trdeg F|K=1$, there are only two
possible cases:
\sn
a) \ $vF=vK$, which implies that $(vF:2vF)=1$,
\n
b) \ $vF=vK\oplus\Z\alpha$ for some $\alpha\in vF\setminus vK$, whence $(vF:2vF)=2$.
\sn
By the Baer--Krull Theorem, in the first case there is no other ordering on $F$ that induces the same place as the
given one. In the second case there is exactly one other ordering that induces the same place.

Now consider the cut that corresponds to the given ordering according to Theorem~\ref{homeo}. In analogy to the case
of a rational function field discussed above, the authors of \cite{KKo} call this cut a \bfind{ball cut} if the
second case holds. The following argument justifies this definition. Pick any element $X\in F\setminus K$. Then
$F|K(X)$ is algebraic, thus $vF/vK(X)$ is a torsion group. This implies that case b) holds for $F$ with the given
ordering if and only if it holds for
$R(X)$ with the restriction of this ordering, as the corresponding natural valuation on $R(X)$ is just the
restriction of the natural valuation on $F$. From this one obtains:

\begin{proposition}
The following are equivalent:
\sn
1) \ the cut corresponding to the given ordering on $F$ is a ball cut,
\n
2) \ for some $X\in F\setminus K$, the cut induced by $X$ in $R$ under the restriction of the ordering to $R(X)$
is a ball cut, 
\n
3) \ for each $X\in F\setminus K$, the cut induced by $X$ in $R$ under the restriction of the ordering to $R(X)$
is a ball cut. 
\end{proposition}

All results above can be obtained in an abstract setting for abstract real curves.
However, once we embed the curve in an affine space we obtain a clearer picture.
Note that every $n$-dimensional affine space ${\mathbb A}^n R$ over $R$ is an ultrametric space with the ultrametric generated
by the natural valuation $v$ of $R$. The ultrametric distance between points $(x_1,...,x_n)$ and $(y_1,...,y_n)$
can be defined and computed as follows:
\[
u((x_1,...,x_n),(y_1,...,y_n)) = \min\{v(x_i - y_i)\} = \frac{1}{2}v\left(\sum (x_i - y_i)^2\right)\>.
\]
Therefore we can consider ultrametric balls in ${\mathbb A}^n R$. We say that an ultrametric ball $B$ in
${\mathbb A}^n R$ cuts a curve $\mathcal C$ if $B\cap \mathcal C \neq \emptyset$ and $({\mathbb A}^n R\setminus B)
\cap \mathcal C \neq \emptyset$. In this case $B$
determines a cut (always more then one) on the curve. In \cite{KKo} it is shown that such a cut is a ball cut,
and the  following theorem is proved:

\begin{theorem}
Every ball cut on a smooth and complete real affine curve in ${\mathbb A}^n R$ is induced by some ultrametric ball.
If the orderings corresponding to two ball cuts induce the same $\R$-place, then there is an ultrametric ball in
${\mathbb A}^n R$ which induces both cuts on the curve.
\end{theorem}

The converse of the second assertion is not true, a counterexample is given in \cite{KKo}. The ball mentioned in
this assertion can well induce more than two cuts on the curve. It is an open question how to determine the pairs
of cuts that induce the same $\R$-place.

\sn
\subsection{Classification of Artin-Schreier defect extensions}                  \label{sectASe}
An \bfind{Artin-Schreier extension} is a field extension $L|K$ of degree $p$ of fields of characteristic $p$
generated by an element $\vartheta$ that satisfies $\vartheta^p-\vartheta \in K$. Such an extension has
nontrivial defect if and only if it is immediate. In this case, the cut $v(\vartheta-K)^+$ taken in the
divisible hull of $vK$
enables us to distinguish two types of Artin-Schreier defect extensions. We call such an extension
\bfind{dependent} if it can be derived by a transformation from a purely inseparable defect extension
of degree $p$, and \bfind{independent} otherwise. In \cite{[Ku4]} the following result is proved:

\begin{theorem}                          \label{charindep}
An Artin-Schreier defect extension is independent if and only if the cut $v(\vartheta-K)^+$ is a
group$^-$-cut.
\end{theorem}

This classification of Artin-Schreier defect extensions is important because work by M.~Temkin (see e.g.\
\cite{[Te]}) and by the author indicates that dependent defect appears to be more harmful to the above cited problems
than independent defect. In the paper \cite{[CP]}, S.~D.~Cutkosky and O.~Piltant give an example of an extension of
valued function fields consisting of a tower of two Artin-Schreier defect extensions where so-called strong
monomialization fails.
As the valuation on this extension is defined by use of so-called generating sequences, it is hard to determine
whether the Artin-Schreier defect extensions are dependent or independent. However, Cutkosky,
L.~Ghezzi and S.~ElHitti show that both of them are dependent (see e.g.\
\cite{[EG]}); this again lends credibility to the hypothesis that dependent defect is the more harmful one.

\pars
Moreover, the classification is an important tool
in the proof of the following theorem in \cite{[Ku4]}:

\begin{theorem}                                \label{chardl}
A valued field of positive characteristic is henselian and defectless if and only if each purely inseparable 
extension is defectless and the field does not allow any proper immediate algebraic extensions.
\end{theorem}

This theorem in turn is used in \cite{[K1]} for the construction of an example showing that a certain natural
axiom system for the elementary theory of $\F_p((t))$ (``henselian defectless valued field of characteristic $p$
with residue field $\F_p$ and value group a $\Z$-group'') is not complete.

\parm
It would be desirable to have a classification of Galois defect extensions of prime degree and an analogue of
Theorem~\ref{chardl} also in the case of valued fields of mixed characteristic (i.e., valued fields of
characteristic $0$ with residue fields of positive characteristic). The obvious problem is to find the
suitable definition of ``dependent'', since there are no purely inseparable extensions. Some guidance can be obtained
from the theory of perfectoid fields, as these allow an exchange of information between the mixed characteristic
case and the case of equal positive characteristic. If we follow this guidance, then all Galois defect extensions
of prime degree of perfectoid fields should be called independent. An indication that this is the
right choice comes from the fact that perfectoid fields are deeply ramified fields in the sense of Section 6.6 of
\cite{GR}; recent work of the author of these notes indicates that when $\vartheta$ is a suitably chosen generator
of any Galois defect extension of prime degree of a deeply ramified field, then as in Theorem~\ref{charindep},
$v(\vartheta-K)^+$ is a group$^-$-cut.

\sn
\subsection{Approximation of elements in henselizations}                  
Complete valued fields of rank 1 (i.e., with archimedean ordered value group) are henselian, but for valuations
$v$ of arbitrary rank this does not hold in general. However, there is a
connection between Hensel's Lemma and completions, but these completions
have to be taken for residue fields of suitable coarsenings of $v$. This
connection was worked out by P.~Ribenboim in \cite{[Ri]} who used \bfind{distinguished
pseudo Cauchy sequences} to characterize the so called \bfind{stepwise
complete} fields; it had been shown by Krull that these fields are henselian.

Take any immediate extension $(L|K,v)$ of valued fields and $z\in L\setminus K$.
We call $z$ \bfind{weakly distinguished over $K$} if $v(z-K)^+$ is a ball$^+$-cut, and we call $z$
\bfind{distinguished over $K$} if it is a group$^+$-cut. The latter name is chosen since distinguished elements
are limits of distinguished pseudo Cauchy sequences.

Now take an arbitrary valued field $(K,v)$ and extend its valuation $v$
to its algebraic closure $\tilde{K}$. Then $\tilde{K}$ contains a unique
henselization $K^h$ with respect to this extension. The following result is proved in \cite{[Ku5]},
answering a question from B.~Teissier. It has recently been reproven by Teissier using methods from
algebraic geometry.
\begin{theorem}                             \label{MT1}
Each element in $K^h\setminus K$ is weakly distinguished over~$K$.
\end{theorem}

Note that if $(K,v)$ is of rank 1, then its henselization lies in its completion and every
element $a\in K^h\setminus K$ is distinguished over $K$ (with $v(a-K)^+=(vK)^+$).

By ``$\alpha>v(a-K)$'' we mean $\alpha>v(a-c)$ for all $c\in K$.
Theorem~\ref{MT1} is used in \cite{[Ku5]} to prove the following result:

\begin{theorem}                             \label{MT2}
Take $z\in \tilde{K}\setminus K$ such that
\[
v(a-z)\;>\;v(a-K)
\]
for some $a\in K^h$. Then $K^h$ and $K(z)$ are not linearly
disjoint over $K$, that is,
\[
[K^h(z):K^h] < [K(z):K]
\]
and in particular, $K(z)|K$ is not purely inseparable.
\end{theorem}

This theorem has a crucial application in \cite{[Ku4]} to the
classification of Artin-Schreier defect extensions which we discussed in the previous section.
The classification was originally obtained in \cite{[Ku1]} under the
additional assumption that the fields in question are henselian.
With the help of Theorem~\ref{MT2} this assumption can be dropped, and
so the classification becomes available for valued function fields.

\mn
%
%
\section{The invariance group}              \label{sectinvgp}
For every cut $\Lambda$ in an ordered abelian group $G$, we define
\[
\cG(\Lambda)\>:=\> \{g\in G\mid \Lambda^L+g=\Lambda^L\}
\]
and call it the \bfind{invariance group} of $\Lambda$; other authors (e.g.\ Ehrlich in \cite{Eh}, following
Kijima and Nishi \cite{[KN]}) call it the \bfind{breadth} of the
cut $\Lambda$. Note that $\Lambda^L+g=\Lambda^L$ is equivalent to $\Lambda^R+g=\Lambda^R$.

\pars
The proof of the following facts is straightforward (see e.g.\ \cite{[Knew]}).
\begin{lemma}                                    \label{invgr}
Take an ordered abelian group $G$ and a Dedekind cut $\Lambda$ in $G$. Then $\cG(\Lambda)$ is a convex subgroup
of $G$, and $G$ is the disjoint union of the three convex subsets $\Lambda^L-\Lambda^R$, $\cG(\Lambda)$ and
$\Lambda^R-\Lambda^L$, with
\[
\Lambda^L-\Lambda^R\><\>\cG(\Lambda)\><\>\Lambda^R-\Lambda^L\>.
\]
\end{lemma}

\begin{corollary}                         \label{invgptriv}
The invariance group of $\Lambda$ is trivial if and only if
\[
\Lambda^R-\Lambda^L=G^{>0}\>.
\]
\end{corollary}

The following theorem is proved in \cite{[Knew]}, but has also been stated (more or less explicitly) by other
authors:
\begin{theorem}                           \label{ball-invg}
A cut $\Lambda$ in an ordered abelian group is a ball cut if and only if it is the upper or lower edge
of a coset of its invariance group, i.e., if there is some $g\in G$ such that $\Lambda=(g+\cG(\Lambda))^+$ or
$\Lambda=(g+\cG(\Lambda))^-$.
\end{theorem}

\begin{remark}                  \label{eigentlich}
\rm R.~Baer (\cite{B}) introduced the notion \bfind{eigentlicher Schnitt}, that is, a Dedekind cut $\Lambda$ in an
ordered abelian group $G$ such that for every positive $g\in G$ there are $a\in \Lambda^L$ and $b\in\Lambda^R$
such that $b-a<g$. By Lemma~\ref{invgr} this condition is equivalent to the invariance group of $G$ being trivial.
Ehrlich (\cite{Eh}) calls them \bfind{Veronese cuts}, and Galanova and Pestov call them \bfind{fundamental cuts}.
The non-principal cuts with trivial invariance group are called \bfind{dense} in Tressl's papers \cite{[T1],[T2]}.
\end{remark}

Several authors, e.g.\ Rolland in \cite{[Ro]} and Wehrung in \cite{[W]}, work with the ball$^+$-cuts
$\cG(\Lambda)^+$ rather than the invariance groups themselves. The set of all of these cuts in an ordered abelian
group $G$ coincides with the set of cuts $H^+$ where $H$ runs through all (proper) convex subgroups of $G$. (Note
that $G$ itself is the invariance group of the two cuts $(G,\emptyset)$ and $(\emptyset,G)$, which are not Dedekind
cuts.)

\pars
If $\Lambda$ is a cut in an ordered field $K$ and is positive (i.e., $0\in\Lambda^L$), then it is also a cut in
the ordered abelian multiplicative group of positive elements of $K$. Its invariance group there is called the
\bfind{multiplicative invariance group of $\Lambda$}, and we denote it by $\cG^\times(\Lambda)$.

\sn
\subsection{Invariance group and pseudo Cauchy sequences}
From what we have said about immediate extensions and pseudo Cauchy sequences in Section~\ref{sectpCs}, ordered
abelian groups or fields that are maximal with respect to their natural valuation contain limits for all pseudo
Cauchy sequences. This is why several authors employ pseudo Cauchy sequences to study and to characterize such
ordered abelian groups or fields. Certain cuts can induce, or be induced by, pseudo Cauchy sequences. For example,
if $(a_\nu)_{\nu<\lambda}$ is a pseudo Cauchy sequences which is also strictly increasing, then it is cofinal in
the lower cut set of the cut $\Lambda=\{a_\nu\mid \nu<\lambda\}^+$, and the following holds:

\begin{theorem}
The invariance group of $\{a_\nu\mid \nu<\lambda\}^+$ is equal to the breadth of the pseudo Cauchy sequence
$(a_\nu)_{\nu<\lambda}$.
\end{theorem}

If the pseudo Cauchy sequence lies in an ordered abelian group $G$, then it induces a Cauchy sequence (i.e., a
pseudo Cauchy sequence with breadth $\{0\}$) in $G/\cG(\Lambda)$.

\sn
\subsection{Ordered fields with maximal natural valuation}          \label{sectmaxf}
%
%
Take an ordered field $K$. We will denote the ordered additive group of $K$ by $K_+\,$.
In \cite{[KN]}, Kijima and Nishi use the invariance group for the following result:

\begin{theorem}                         \label{KNthm}
The following assertions are equivalent:
\sn
1) \ the natural valuation of $K$ is maximal and its residue field is $\R$,
\sn
2) \ for each cut $\Lambda=(\Lambda^L,\Lambda^R)$ in $K$, the induced cut $(\Lambda^L/\cG(\Lambda),
\Lambda^R/\cG(\Lambda))$ in the ordered abelian group $K_+/\cG(\Lambda)$ is principal.
\end{theorem}
Here, the induced cut is $(\{a/\cG(\Lambda)\mid a\in\Lambda^L\}\,,\,\{b/\cG(\Lambda)\mid b\in\Lambda^R\})$; note
that the two sets are disjoint by the defining property of $\cG(\Lambda)$.

This theorem also holds for any ordered abelian group $G$ in place of the ordered field $K$ if we replace
``its residue field is $\R$'' by ``all of its archimedean components are isomorphic to $\R$''; see the next section.

\parm
In his thesis \cite{[H]} H\"uper considers ordered fields with arbitrary compatible valuations (i.e., valuations
whose valuation ring is convex, or equivalently, contains the valuation ring of the natural valuation). We will
cite one of his main results; in its formulation he uses a notion that is derived from Baer's ``eigentlicher
Schnitt'' (see Remark~\ref{eigentlich}) without explicitly using invariance groups. But using them as follows puts
the result in a wider context:

\begin{theorem}                         \label{HThm}
Take an ordered field $K$ with a compatible valuation $v$. Then the following assertions are equivalent:
\sn
1) \ the valuation $v$ is maximal, 
\sn
2) \ if $H$ is a $\cal O_v$-submodule of $K$ not contained in a larger $\cal O_v$-submodule $H'$ such that
there is no $\cal O_v$-submodule properly between $H'$ and $H$, and if $\Lambda$ is a cut such that
$\Lambda/H:=(\Lambda^L/H,\Lambda^R/H)$ is a Dedekind cut in $K_+/H$ with trivial invariance group, then
$\Lambda/H$ is principal.
\end{theorem}

Let us evaluate this theorem for the case of $v$ being the natural valuation. In this case, condition 2) can be reformulated
as follows:
\sn
{\it 2') \ if $H$ is a convex subgroup of $K$ which is not contained in a larger convex subgroup $H'$ such that
$H'/H$ is archimedean ordered, and if $\Lambda$ is a cut such that $\Lambda/H:=(\Lambda^L/H,\Lambda^R/H)$ is a
Dedekind cut in $K_+/H$ with trivial invariance group, then $\Lambda/H$ is principal.}
\sn
Condition 2') can be further reformulated and thereby simplified by using the following two facts:

\begin{lemma}                      \label{modH}
Take a Dedekind cut $\Lambda$ in an ordered abelian group $G$ and a proper convex subgroup $H$ of $G$. Then the
following assertions hold.
\sn
a) \ $\Lambda/H$ is a Dedekind cut in $G/H$ if and only if $H\subseteq \cG(\Lambda)$.
\sn
b) \ $\cG(\Lambda/H)=\{0\}$ if and only if $\cG(\Lambda)\subseteq H$.
\end{lemma}

In view of these facts, condition 2) is equivalent to:
\sn
{\it 2'') \ if $H$ is a convex subgroup of $K_+$ which is not contained in a larger convex subgroup $H'$ such that
$H'/H$ is archimedean ordered, and if $\Lambda$ is a cut with invariance group $H$, then $\Lambda/H$ is principal.}

\pars
What is the role of the assumption on $H$ in conditions 2') and 2'')? Well, if $H'$ is a larger convex subgroup of
$K_+$ such that $H'/H$ is archimedean ordered, then $H'/H$ is an archimedean component of $K_+$ and therefore
isomorphic to the additive group of $Kv$. As the theorem does not assume
that the latter is equal to $\R$, $H'/H$ may have a non-principal Dedekind cut, which then gives rise to a
non-principal Dedekind cut of $G/H$. So if we are only interested in maximality, then we have to take this case
into account. However, if we assume in addition to condition 1) that $Kv$ is equal to $\R$, then all
Dedekind cuts in archimedean components are principal, and we can drop the assumption on $H$. This shows that
Theorem~\ref{KNthm} is a consequence of Theorem~\ref{HThm}.

\pars
\begin{remark}
\rm In \cite{[KNS]}, the authors state that ``the notion of maximal ordered fields was first introduced'' in earlier
papers of theirs, the earliest published in 1987 by Kijimi and Nishi. This statement is correct only as far as it concerns results published in journals, as the work of H\"uper shows.
\end{remark}

\sn
\subsection{Archimedean complete ordered abelian groups}          \label{sectmaxg}
An ordered abelian group $(G,<)$ is called \bfind{archimedean complete} (a notion introduced by H.~Hahn) if every
proper ordered abelian group extension $(G',<)$ of $(G,<)$ introduces new archimedean classes, or in other words,
the natural valuation on $(G',<)$ has a larger value set than on $G$. Hence the archimedean complete
ordered abelian groups are precisely the ordered abelian
groups that are maximal w.r.t.\ their natural valuation and whose archimedean components are as large as possible,
that is, isomorphic to the additive group of real numbers. Hahn shows in \cite{Ha} that archimedean
complete ordered abelian groups are precisely the ones that admit an order preserving isomorphism onto a so-called
Hahn product with all of its archimedean components equal to the additive group of real numbers. (Hahn products
are the analogues for ordered abelian groups of the power series fields.)

Archimedean complete ordered abelian groups $G$ are characterized in the paper \cite{CG} by L.~W.~Cohen and C.~Goffman
as follows:
\begin{theorem}
An ordered abelian group $(G,<)$ is archimedean complete if and only if for every proper convex subgroup $H$,
the ordered factor group $G/H$ is dense and every cut in $G/H$ with trivial invariance group is principal.
\end{theorem}

Ehrlich revisits this topic in \cite{Eh}. Relying on the Theorem of Cohen and Goffman, Ehrlich proves:
\begin{theorem}
An ordered abelian group $(G,<)$ is archimedean complete if and only if for every cut $\Lambda$, the induced cut
in $G/\cG(\Lambda)$ is principal, but not a jump.
\end{theorem}
Since ordered fields admit no jumps, this theorem can be seen as an analogue of Theorem~\ref{KNthm}.
Ehrlich shows that the induced cut has trivial invariance group; this is a special case of part b) of
Lemma~\ref{modH}.

\parm
We recommend Ehrlich's paper  \cite{Eh} for interesting historical remarks and a detailed list of references.

\sn
\subsection{Model theory of ordered fields with cuts}                  \label{sectTres}
In \cite{[T],[T1],[T2]}, Tressl studies the model theory of real closed fields with a fixed cut.
Given a model $M$ of an o-minimal extension $T$ of the theory of real closed fields in a language $\mathcal L$,
he determines the model theoretic properties of $M$ in the language ${\mathcal L}({\mathcal D})$ where $\mathcal D$
is a predicate for the left cut set $\Lambda^L$ of a fixed cut $\Lambda$ in $M$. If $(M_1,\Lambda^L_1)$ and
$(M_2,\Lambda^L_2)$ are two structures obtained in this way, conditions are found for $(M_1,\Lambda^L_1)$
and $(M_2,\Lambda^L_2)$ to be elementarily equivalent in the language ${\mathcal L}({\mathcal D})$ enhanced by
parameters from a common elementary substructure of $M$ and $M'$. The main result is rather technical in nature,
but for special classes of cuts, the situation is much easier. To illustrate this, the following theorem is taken
from \cite{[T2]}:

\begin{theorem}
Let $A\prec M_1,M_2$ be models of $T$ and let $\Lambda_1,\Lambda_2$ be non-principal cuts in $M_1,M_2$
respectively, with trivial invariance groups. Then $(M_1, \Lambda^L_1) \equiv_A (M_2, \Lambda^L_2)$ if and only if
the restrictions of $\Lambda_1$ and $\Lambda_2$ to $A$ coincide.
\end{theorem}

\mn
%
%
\section{The invariance valuation ring}              \label{sectinvvr}
The \bfind{invariance valuation ring of a cut} $\Lambda$ in an ordered field $K$ is defined as
\[
{\cal O}(\Lambda)\>:=\>
\{b\in K\mid b\,{\cal G}(\Lambda)\subseteq {\cal G}(\Lambda)\}\>.
\]
We denote its maximal ideal $\{b\in K\mid b\,{\cal G}(\Lambda)\subsetuneq {\cal G}(\Lambda)\}$ by
${\cal M} (\Lambda)$.

According to Lemma~\ref{invgr} ${\cal G}(\Lambda)$ is a convex subgroup of the ordered additive group $K_+$
of $K$, and we have already noted in
the Introduction that every convex subgroup of $K_+$ is an ${\cal O}_v$-module, where $v$ denotes the natural
valuation. In this way, the above definition becomes a special case of the following.

Take any valued field $K$ with valuation ring ${\cal O}_v$ and an ${\cal O}_v$-module $M$ in $K$.
The \bfind{invariance valuation ring of an ${\cal O}_v$-module} $M$ in $K$ is defined as
\[
{\cal O}(M)\>:=\>
\{b\in K\mid bM\subseteq M\}\>.
\]

The relation between multiplicative invariance group and invariance valuation ring is determined in \cite{[Knew]}. Also Tressl and G\"uldenberg obtain results on this topic.

\sn
\subsection{Projecting cuts into residue fields}              \label{sectproj}
Take a convex valuation ring ${\cal O}$ of an ordered field $K$, with maximal ideal
${\cal M}$. Its residue field ${\cal O}/{\cal M}$ is again an ordered
field, with the ordering induced through the residue map. We will say that
the cut $\Lambda$ \bfind{can be projected into the residue field
${\cal O}/{\cal M}$} if there are elements $a,c\in K$ such that
$c>0$ and $c\Lambda+a$ induces a Dedekind cut
\begin{equation}                            \label{rescut}
\left(\,((c\Lambda^L+a)\cap {\cal O})/{\cal M}\,,\,
((c\Lambda^R+a)\cap {\cal O})/{\cal M}\,\right)
\end{equation}
in ${\cal O}/{\cal M}$ via the residue map.

The following theorem shows for which convex valuation rings $\cal O$ a cut can be projected into the associated
residue field. For a proof, see \cite{[Knew]}.

\begin{theorem}                             \label{projrf}
1) \ Take any convex valuation ring ${\cal O}$ of $(K,<)$. If ${\cal O}
(\Lambda)\subsetuneq {\cal O}$, then the cut $\Lambda$ can be
projected into the residue field ${\cal O}/ {\cal M}$. If ${\cal O}
\subsetuneq {\cal O}(\Lambda)$, then it cannot be projected into
${\cal O}/{\cal M}$.\n
2) \ The cut $\Lambda$ can be projected into ${\cal O}(\Lambda)/
{\cal M}(\Lambda)$ if and only if $(v{\cal G}(\Lambda))^-$
is a ball$^{\,+}$-cut.
\end{theorem}

\sn
\subsection{Definable valuation rings in ordered fields}              \label{sectdefval}
It is obvious that if a cut (that is, its lower cut set $\Lambda^L$) is definable in some extension of the language
of ordered rings, then so is $\cG(\Lambda)$. It then follows that also the invariance valuation ring
is definable.

This observation is put to work in \cite{JSW}, where the following is proved:

\begin{proposition}
Take an ordered field $K$ with real closure $R$. If $K$ is not dense in $R$, then $K$ admits a nontrivial
valuation definable in the language of ordered rings.
\end{proposition}
\sn
The idea of the proof is as follows. If $K$ is not dense in $R$, then there is an element $r\in R\setminus K$ and a
positive element $a\in R$ such that $|r-c|>a$ for all $c\in K$. Since $R|K$ is algebraic, the set $K^{>0}$ of
positive elements in $K$ is coinitial in $R^{>0}$, so we can choose $a\in K$.
If we set $\Lambda^L=\{c\in K\mid c<r\}$, then we obtain a cut
$\Lambda$ such that $\Lambda^R-\Lambda^L\subsetuneq K^{>0}$. By Corollary~\ref{invgptriv}, its invariance group is
thus nontrivial. This implies that the invariance valuation ring is not all of $K$, so the associated valuation is
nontrivial. Since $r$ lies in a real closure of $K$, the set $\Lambda^L=\{c\in K\mid c<r\}$ is definable, and by
what we said above, so are the invariance valuation ring and thus also the associated valuation.

The above arguments also prove the following general principle:\n
{\it If the lower cut set of some Dedekind cut
with nontrivial invariance group in an ordered field is definable, then the field admits a nontrivial definable
valuation ring.}

Similarly, if an ${\cal O}_v$-module $M$ in a valued field $(K,v)$ is definable, then so is its
invariance valuation ring ${\cal O}(M)$. If $K\ne M\ne\{0\}$, then ${\cal O}(M)$ is a nontrivial valuation ring.
This yields the following general principle: \n
{\it If a proper nontrivial ${\cal O}_v$-module in a valued field
$(K,v)$ is definable, then the field admits a nontrivial definable valuation ring containing~${\cal O}_v\,$.}
\sn
Note that if a cut $\Lambda$ is definable in the value group $vK$ in a suitable language of valued fields, then the
${\cal O}_v$-module $\{a\in K\mid va\in \Lambda^R\}$ is also definable, and it is proper and nontrivial if and only
if the cut is a Dedekind cut.

\mn
%
%
\section{Cut cofinalities}                   \label{sectcof}
Recall that by the \bfind{cofinality} of the cut $\Lambda$ we mean the pair
$(\kappa,\lambda)$ where $\kappa$ is the cofinality of $\Lambda^L$, and $\lambda$ is the
coinitiality of $\Lambda^R$.

\sn
\subsection{Ordered abelian groups and fields that are $\eta_\alpha$-sets}       \label{sectea}
Take any ordinal $\alpha$. An $\eta_\alpha$-set is an ordered set $S$ such that for any two subsets $A\subseteq S$
and $B\subseteq S$ of cardinality less than $\aleph_\alpha$ with $A<B$, there is some $s\in S$ such that $A<s<B$.
This is equivalent to saying that $S$ does not admit any cuts of cofinality $(\kappa,\lambda)$ where both
$\kappa$ and $\lambda$ are smaller than $\aleph_\alpha\,$. In \cite{A}, Alling proves:

\begin{theorem}
a) An ordered abelian group is an $\eta_\alpha$-set if and only if it is $\alpha$-maximal, its value set w.r.t.\
the natural valuation is an $\eta_\alpha$-set, and all of its archimedean components are isomorphic to $\Z$ or $\R$.
\sn
b) An ordered field is an $\eta_\alpha$-set if and only if it is $\alpha$-maximal, its value group w.r.t.\ the
natural valuation is an $\eta_\alpha$-set, and its residue field is $\R$.
\end{theorem}

Every $\aleph_\alpha$-saturated ordered abelian group or field is an $\eta_\alpha$-set. For the converse, the
reader may note that divisible ordered abelian groups and real closed fields are o-minimal. This implies that for
them, the property of being an $\eta_\alpha$-set is equivalent to that of being $\aleph_\alpha$-saturated. For
results on the variety of $\eta_\alpha$ ordered abelian groups or fields, for fixed $\alpha$, see \cite{AK}.

\pars
Ball cuts and invariance groups do not appear explicitly in \cite{A} or \cite{AK}. But Rolland draws a connection
in \cite{[Ro]}. He states that an ordered abelian group is an $\eta_\alpha$-set if and only if its value set
is an $\eta_\alpha$-set and for every Dedekind cut $\Lambda$ with nontrivial invariance group $\cG(\Lambda$, the
coinitiality of the upper cut set of $\cG(\Lambda)^+$ is not less than $\aleph_\alpha$.

\sn
\subsection{Cuts in ordered power series fields}                  \label{sectgp}
The papers \cite{G1,G2,G3,GP} of Galanova and Pestov are devoted to the study of cuts in power series fields and
in restricted power series fields (in the latter, the cardinalities of the supports of the power series are
bounded by a given cardinal number). We cite three theorems from \cite{GP}. The cardinality of a set $S$ is
denoted by $|S|$, and $|S|^+$ denotes its successor cardinal.

\begin{theorem}
Take any ordered abelian group $G$. Then all cuts in the power series field $\R((G))$ are ball cuts.
\end{theorem}
The proof of this theorem in \cite{GP} is long and technical. Let us give the sketch of a shorter and more
conceptual proof. We write $K=\R((G))$. Every cut $\Lambda$ in $K$ is realized in some ordered field extension
$L$ of $K$
(for instance, if $L$ is a $|K|^+$-saturated elementary extension of the ordered field $K$). As a power series
field, $K$ is maximal w.r.t.\ its natural valuation. Extend $v$ to the natural valuation of $L$. Then it
follows that for every $x\in L\setminus K$ there is some $a\in K$ such that $v(x-a)=\max\{v(x-c)\mid c\in K\}$
since otherwise, $x$ would be a limit of some pseudo Cauchy sequence in $K$ without a limit in $K$, contradicting
the fact that $K$ is maximal. The value $v(x-a)$ can only be maximal if either $v(x-a)\notin vK$ or there is
$d\in K$ such that $vd=v(x-a)$ and $d^{-1}
(x-a)v\notin Kv$. But the latter cannot be the case: since $Kv=\R$ and $v$ on $L$ is a natural valuation, we
must have that $Lv=Kv$. Hence $\gamma:=v(x-a)\notin vK$. We leave it as an exercise to the reader to show that
$\Lambda=(a+\{b\in K\mid vb >\gamma\})^+$ or $\Lambda=(a+\{b\in K\mid vb >\gamma\})^-$.

\begin{theorem}
Take any ordered abelian group $G$ and a cardinal number $\kappa$ such that $\aleph_0<\kappa<|G|$.
Denote by $\R((G,\kappa))$ the subfield of $\R((G))$ consisting of all power series with support of cardinality
less than $\kappa$. Take a non-ball cut in $\R((G,\kappa))$ of cofinality $(\lambda,\lambda)$. Then $\lambda$ is
equal to the cofinality of $\kappa$. In particular, if $\kappa$ is regular, then $\lambda=\kappa$.
\end{theorem}

\begin{theorem}
Take a non-principal cut $\Lambda$ in some ordered field $K$ with trivial invariance group, and let $(\kappa,
\kappa)$ be its cofinality. Then $\kappa$ is equal to the cofinality of $K$.
\end{theorem}

\parm
In \cite{[Ro]}, Rolland states the existence of power series fields that admit cuts with preassigned cofinalities
$(\kappa_i,\lambda_i)$, $i\in I$, where the $\kappa_i$ and $\lambda_i$ are infinite regular cardinals. The proof
he gives is insufficient, but the result also follows from the work we will discuss in the next section. He uses
it to show the existence and to (partially) characterize the ordered fields which admit a closed bounded interval
and a continuous function which is unbounded on this interval.

\sn
\subsection{Symmetrically complete ordered abelian groups and fields}                  \label{sectsc}
A Dede\-kind cut with cofinality $(\kappa,\lambda)$ is called \bfind{symmetric} if $\kappa=\lambda$, and
\bfind{asymmetric} otherwise. Note that the notion of symmetry used by Galanova and Pestov is quite different
from the one defined here. However, Pestov states in \cite{Pe}, without proof or reference) that if a cut in an
ordered field is symmetric in their sense (i.e., it is not a ball cut), then it is also symmetric in the sense of
the above definition. A proof is given by Galanova in \cite{G1}. In \cite{[Ro]}, Rolland states the same
in full generality for non-ball cuts in ordered abelian groups. The statement is correct, but his proof appears to
have a serious gap.

\pars
A linearly ordered set $(I,<)$ is called \bfind{symmetrically complete} if every symmetric cut in $I$ has
cofinality $(1,1)$, i.e., is a jump. In dense linear orderings (and hence in ordered fields) there are no jumps.
Consequently, a dense linear ordering is symmetrically complete if and only if all of its cuts are asymmetric.

For example, $\Z$ and $\R$ are symmetrically complete, but $\Q$ is not.
In $\Z$ and $\R$, every Dedekind cut is principal; in $\Z$ all of them have
cofinality $(1,1)$, and in $\R$ they have cofinalities $(1,\aleph_0)$
and $(\aleph_0,1)$. In contrast, in $\Q$ the Dedekind cuts have cofinalities
$(1,\aleph_0)$, $(\aleph_0,1)$ and $(\aleph_0,\aleph_0)$.

\pars
In \cite{KKS} it is shown that a symmetrically complete ordered abelian group is spherically complete w.r.t.\ its
natural valuation and hence a Hahn product, with all of its archimedean components isomorphic to $\R$. Similarly, a
symmetrically complete ordered field is spherically complete w.r.t.\ its natural valuation and hence a power
series field, with residue field $\R$. For Hahn products with all of its archimedean components isomorphic to
$\R$ the set of all cut cofinalities is computed from the set of all cut cofinalities of the value set of its
natural valuation, and a similar computation is done for power series fields with residue field $\R$.
Based on this computation, a full characterization of symmetrically complete ordered abelian groups and fields
is obtained. We will cite a selection of the main results.

We call
an ordered set \bfind{strongly symmetrically complete} if it is symmetrically complete and does not have any cuts
with cofinalities $(1,\aleph_0)$ or $(\aleph_0,1)$.

\begin{theorem}
A non-archimedean ordered field is symmetrically complete if and only if it is
spherically complete w.r.t.\ its natural valuation, has residue
field $\,\R$ and a dense strongly symmetrically complete value group.

\pars
A nontrivial densely ordered abelian group is symmetrically
complete if and only if it is spherically complete w.r.t.\ its natural
valuation $v$, has a dense strongly symmetrically complete value set, and all archimedean components are isomorphic
to $\R$. It is strongly symmetrically complete if and only if in addition,
the value set has uncountable cofinality.
\end{theorem}
\n
In particular, it follows that symmetrically complete ordered abelian groups are divisible and symmetrically complete ordered fields are real closed.

\pars
Further, it is shown in  \cite{KKS} that every ordered set can be extended to a dense strongly symmetrically
complete ordered set with uncountable cofinality. The reader may note that this result is not explicitly stated in
Hausdorff's work. The authors of \cite{KKS} also tried to give a construction that is as short as possible. It
turns out that the constructed orderings are themselves lexicographic, as are the orderings on Hahn products and
power series fields. Such orderings deserve to be studied in more detail.

Using the above results and the fact that every ordered set is the natural value set of
some Hahn product with all components isomorphic to $\R$, and every ordered abelian group is the natural value
group of some power series field with residue field $\,\R$, the following result of \cite{[S]} is reproved:

\begin{theorem}
Every ordered field can be embedded in a symmetrically complete ordered field. Every ordered abelian group can be
embedded in a symmetrically complete ordered abelian group.
\end{theorem}


\end{document}